\RequirePackage{fix-cm}
\documentclass[smallextended]{svjour3}       % onecolumn (second format)
\smartqed  % flush right qed marks, e.g. at end of proof

\usepackage{array}
\usepackage{amssymb}
\usepackage{booktabs}
\usepackage{subfigure}
\usepackage{graphicx}
\begin{document}

\title{A new approach for solving nonlinear Thomas-Fermi equation based on fractional order of rational Bessel functions}


\titlerunning{Solving nonlinear Thomas-Fermi equation}        % if too long for running head

\author{K. Parand$^{\dagger,\ast}$ \and A. Ghaderi$^{\ast}$ \and M. Delkhosh$^{\ast}$ \and H. Yousefi$^{\ast}$}


\institute{$\ast$\at Department of Computer Sciences, Shahid Beheshti University, G.C., Tehran, Iran.\\ \at {$\dagger$} Department of Cognitive Modelling, Institute for Cognitive and Brain Sciences, Shahid Beheshti University, Tehran, Iran, Email: k\_parand@sbu.ac.ir.
}

\date{Received: date / Accepted: date}

\maketitle

\begin{abstract}
In this paper,  the fractional order of rational Bessel functions collocation method (FRBC) to solve Thomas-Fermi equation which is defined in the semi-infinite domain and has singularity at $x = 0$ and its boundary condition occurs at infinity, have been introduced. We solve the problem on semi-infinite domain without any domain truncation or transformation of the domain of the problem to a finite domain. This approach at first, obtains a sequence of linear differential equations by using the quasilinearization method (QLM), then at each iteration solves it by FRBC method. To illustrate the reliability of this work, we compare the numerical results of the present method with some well-known results in other to show that the new method is accurate, efficient and applicable.

\end{abstract}
\section{Introduction}\label{sec1}

Many problems arising in mathematical, fluid dynamics, quantum mechanics, astrophysics, physics and engineering are occurred infinite or semi-infinite domains. In this section, will be expressed about some of approaches for solving problems which are defined in unbounded domains and  a brief history of Thomas-Fermi equation that is defined on semi-infinite domain.

\subsection{Solving problems  over unbounded domains}
Recently various approaches have been successfully proposed to solve problems which are occurred  unbounded domains. Such as numerical, analytical and semi-analytical methods. 

Different numerical methods have been introduced for problems  which is defined in the semi-infinite domain such as finite difference method (FDM) \cite{refdd1,refdd2}, finite element method (FEM) \cite{refdd2,refdd3}, meshfree methods \cite{refdd4,refdd5}, and spectral methods  \cite{refdd6}.

The study of analytical and semi-analytical solutions of differential equations (DEs) plays an important role in mathematical physics, engineering and the other sciences. In the past several decades, various methods for obtaining solutions of DEs have been presented, such as Adomian decomposition method \cite{refdd19,refdd20}, Homotopy perturbation method \cite{refdd21}, Variational iteration method \cite{refdd22}, Exp-function method \cite{refdd23} and so on.

The spectral approximations for DEs on finite domains have achieved great success and popularity in recent years, but spectral approximations for DEs on infinite domains have only received limited attention. Several spectral methods for treating infinite/semi-infinite domain problems have been utilized by different researchers: (1) There is an effective approach to solve these problems by applying  the basis of Sinc, Hermite and rational Christov functions that are orthogonal on the interval $(-\infty,\infty)$ \cite{refdd7,refdd8,refdd9,refdd10}  and basis of Laguerre polynomials that  are orthogonal on the interval $[0,\infty)$ \cite{refdd11,refdd12}. (2) Another approach for solving such problems which is based on rational approximations. This method transfers polynomials on interval $[\alpha,\beta]$ to functions on interval $[0,\infty)$ by using the algebraic mapping $x\to \frac {\beta x+\alpha  L}{x + L}$ that $L>0$ is a scaling/stretching factor \cite{refdd13,refdd14}. The Jacobi polynomials are a class of classical orthogonal polynomials, also the Gegenbauer polynomials, the Legendre and Chebyshev polynomials, are special cases of these polynomials which have been used in several literatures for solving some problems. In some papers have been provided the collocation method for natural convection heat transfer equations embedded in porous medium, nonlinear differential equation and nonlinear integro-differential equation based on rational Gegenbauer, Legendre, Chebyshev functions \cite{refdd8,refdd11,refdd24,refdd25}.  Doha et al have presented Jacobi rational-Gauss collocation method that based on Jacobi rational functions and Gauss quadrature integration to solve nonlinear Lane-Emden equation \cite{refdd16}. Isik et al have used Bernstein polynomials to solve high order initial and boundary values problems. Their approximate solution has better convergence rate than the one found by using collocation method \cite{refdd15}. (3) Guo \cite{Guo1,Guo2} applied a method that proceeds by mapping the original problem in an unbounded domain to a problem in a bounded domain, and then using suitable Jacobi polynomials such as Gegenbauer polynomials to approximate the resulting problems. (4) A further approach consists of replacing the infinite domain with $[-K,K]$ and the semi-infinite domain with $[0,K]$ by choosing $K$ sufficiently large. This method is named domain truncation \cite{reff001,refdd17,refdd18}.

In this investigation, we attempt to introduce a spectral method based on the fractional order of rational Bessel functions (FRB) to solve Thomas-Fermi on semi-infinite domian.
\subsection{The Thomas-Fermi equation}
One of the most important nonlinear ordinary differential equations that occurs in semi-infinite interval is Thomas-Fermi equation, as following \cite{refd00,refd01,refd02}:
\begin{equation}\label{eqqq01}
\frac{d^2y(x)}{dx^2}-\frac{1}{\sqrt{x}}y^{\frac{3}{2}}(x)=0,~~~~~~x\in[0,\infty),
\end{equation}
where boundary conditions for this equation are as follows:
\begin{equation}\label{eqqq02}
y(0)=1,~~~~~\lim_{x\rightarrow\infty}y(x)=0.
\end{equation}

The Thomas-Fermi equation appears in problem of determining the effective nuclear charge in heavy atoms, and because of its importance to theoretical physics, computing its solutions has attracted the attention of the Nobel laureates John Slater (chemistry) \cite{refd03} and Richard Feynman (physics) \cite{refd04} and of course Enrico Fermi \cite{refd05}.

One measure of the rapidity of convergence of the procedure is provided by calculation of the value of the initial slope $y'(0)$ of the Thomas-Fermi potential \cite{refd06}. The problem is useful for calculating form-factors and for obtaining effective potentials which can be used as initial trial potentials in self-consistent field calculations. The initial slope $y'(0)$ is difficult to compute by any means, and plays an important role in determining many physical properties of the Thomas-Fermi atom. It determines the energy of a neutral atom in the Thomas-Fermi approximation:
\begin{equation}\label{eqqq03}
E=\frac{6}{7}\left(\frac{4\pi}{3}\right)^{\frac{2}{3}}Z^{\frac{7}{3}}y'(0),
\end{equation}
where $Z$ is the nuclear charge.

For these reasons, the problem has been studied by many researchers and by the different techniques have been solved, that a number of them are as follows:

Baker in 1930 \cite{refd07} has studied the singularity of this equation, and an analytical solution is calculated as follows:
\begin{eqnarray}
y(x)=1-Bx+\frac{4}{3}x^{\frac{3}{2}}-\frac{2}{5}Bx^{\frac{5}{2}}+\frac{1}{3}x^3+\frac{3}{70}B^2x^{\frac{7}{2}}-\frac{2}{15}Bx^4+\cdots, \nonumber
\end{eqnarray}
where $-B$ is the value of the first derivative at the origin that has calculated $y'(0)=-B=-1.588558$.

Esposito in 2002 \cite{refd14} has reported an original method, due to Majorana, that leads to a semi-analytical series solution of the Thomas-Fermi equation with appropriate boundary conditions in terms of only one quadrature, and proved that the series expansion is uniformly convergent in the interval $[0,1]$, and has calculated $y'(0)=-1.588$.

Liao in 2003 \cite{refd16} has employed the Homotopy analysis method and gave an explicit analytic solution of the Thomas-Fermi equation and the related recurrence formula of constant coefficients. The corresponding $m$th-order approximation is
\begin{eqnarray}
y(x)=\sum_{k=0}^{m}\sum_{n=1}^{4k+1}\alpha_{k,n}(1+x)^{-n}, \nonumber
\end{eqnarray}
where $\alpha_{k,n}$ is defined in Eq. 26 the Ref. \cite{refd16}. He has calculated $y'(0)=-1.58712$.

Kobayashi et al in 1955 \cite{refd31} have examined the asymptotic solution of obtained by Coulson and March \cite{refd28}, and have improved their solution:
\begin{eqnarray}
y(x)=\frac{144}{x^3}\left(1-z+0.6256974977z^2-0.3133861150z^3+0.1373912767z^4-\cdots \right), \nonumber
\end{eqnarray}
where $z=\frac{F}{x^c}$, and $F=13.27097391$ and $c=0.7720018726$, and have calculated $y'(0)=-1.588070972$.

Adomian in \cite{refh1} for solving Thomas-Fermi equation proposed a standard decomposition method that was introduced. Briefly as follows:
\begin{eqnarray}
\nonumber y(x)=c_{1}+c_{2}x+L^{-1}x^{\frac{-1}{2}}\sum\limits_{n=0}^{\infty}A_{n}
\end{eqnarray}
where $L^{-1}$ denotes a two-fold integration, $A_n$ denotes the Adomian polynomials generated for $y^{\frac{3}{2}}$, and $c_1$, $c_2 $ are constants of integration. The Adomian decomposition method employs the recursive relation.

Marinca and Herianu in 2011 \cite{refh4} have used a new method  to find analytical approximate solution to Thomas-Fermi equation and called it the Optimal Parametric Iteration Method (OPIM) that this new iteration approach provides us with a convenient way to optimally control the convergence of the approximate solution. This new iteration approach containing a new iteration scheme involves the presence of a finite number of initially unknown  parameters, which are optimally determined. In this way, the approximative initial slope is $y'(0)=-1.5880659888022421$.

Zhu et al in 2012 \cite{refh5}  approximated the original Thomas-Fermi equation by a nonlinear free boundary value problem (FBVP) and applied an iterative method to solve the FBVP. They transformed the FBVP to a nonlinear singular BVP defined on [0, 1] by a change of variables and also employed an adaptive finite element method based on moving mesh to obtain the best approximate solution at each iteration. Best approximation obtained by this method is $y'(0)=-1.58794357$.

A simple and more precise solution to the Thomas-Fermi equation is obtained by making use of the famous Ritz Variational method. Oulne in 2011 \cite{refh7} used a new simple Variational solution of the Thomas-Fermi equation which reproduces the numerical solution accurately in wide range with a correct asymptotic behavior at long distances from the origin and which allows us to calculate with exactness the initial slope. The proposed solution will be developed in power series which have the same form as series solutions that have been obtained previously by Baker \cite{refd07}. In this method, the approximative initial slope is $y'(0)=-1.588071034$.

Abbasbandy and Bervillier in 2011 \cite{reff10} have compared three methods based respectively on Taylor (Maclaurin) series, Pad\'{e} approximates and conformal mappings. $y'(0) = -1.5880710226113753127189 \pm 7 * 10^{-22}$ has been obtained by using the Pad\'{e}-Hankel method.

Boyd in 2013 \cite{reff01} applied collocation method  based on the rational Chebyshev functions on semi-infinite intervals $TL_n(y; L)$ which $L$ is a user-choosable numerical. Boyd employed Newton-Kantorovich iteration to reduce nonlinear differential equation to a sequence of linear differential equations, and he has calculated$$y'(0) = -1.5880710226113753127186845$$with $L=64$ and $m = 600$ collocation points.

MacLeod in 1992 \cite{reff03} has used  two differing approximations on Chebyshev polynomial according to behavior Thomas-Fermi function, one for small $x<40$, one for large $x$. In this method, the approximative initial slope is $y'(0)=-1.5880710226$.

Tavassoli Kajan et al in 2013 \cite{reff04} proposed the rational third-kind Chebyshev pseudospectral method that is indeed a combination of Tau and collocation methods. They reduces the solution of this problem to the solution of a system of algebraic equations. In this method, the approximative initial slope is $y'(0)=-1.588071022611374$. 

Parand et al in \cite{reff05,reff07,reff09} have proposed collocation method on rational Chebyshev,  Hermite polynomials and Sinc functions to solve Thomas-Fermi on semi-infinite interval without truncating it to a finite domain. These methods reduce the solution of this problem to solution of a system of algebraic equations.

Jovanovic et al in 2014 \cite{reff06} solved the Thomas-Fermi equation by applying a spectral method using an exponential basis set in a semi-infinite domain. The goal of the spectral method approach is to find the values of coefficients $a_i$ that best satisfy the following
equation:
\begin{eqnarray}
y(x)=\sum\limits_{i=1}^{N} a_iR_i ,~~~R_i=e^{-\beta_i x}, \nonumber
\end{eqnarray}
where values of $R_i$  are selected in an intuitive way to cover all the possible decay rate. They have reported  detailed about the convergence rate of the initial slope $y'(0)$ for an exponential basis set.

Liu and Zhu in 2015 \cite{reff08} have proposed an iterative method based on Laguerre pseudospectral approximation which the solution of the Thomas-Fermi equation as sum of two parts due to its singularity at the origin. One "singular" part is a power series expansion. The other "smooth" part satisfies a nonlinear two-point boundary value problem. In this method, the approximative initial slope is $y'(0)=-1.588072$. 

Yao in 2008 \cite{refm12} has solved the Thomas-Fermi equation with a kind of analytic technique, named Homotopy analysis method and his answer is $y'(0)=-1.588004950$.

Amore et al in 2014 \cite{refm16} have obtained highly accurate solutions to the Thomas-Fermi equations for atoms and atoms in very strong magnetic fields. And they apply the Pad\'{e}-Hankel method, numerical integration, power series with Pad\'{e} and Hermite-Pad\'{e} approximates and Chebyshev polynomials. They solved Thomas-Fermi for different $x$ and obtain answer for $y(x)$ and $y'(x)$. Their best answer is  $y'(0)=-1.588071022611375312718684509$.

Fernandez in 2011 \cite{refm17} has showed that a simple and straightforward rational approximation to the Thomas-Fermi equation provides the slope at origin with unprecedented accuracy and that Pad\'{e} approximates of relatively low order are far more accurate than more elaborate approaches proposed recently by other authors. He calculated $y(x)$ for different $x$ and compare their method with Chebyshev and numerical method, and calculated  $y'(0)=-1.588071022611375313$.

Epele et al in 1999 \cite{refm18} have used Pad\'{e} approximate approach to solving Thomas-Fermi equation. They have calculated $y'(0)=-1.5881$. 

Khan and Xu in 2007 \cite{refm19} have used an analytic technique, namely the Homotopy analysis method (HAM). Their best answer for $y'(0)$ was $-1.586494973$ when they selected [30,30] for Homotopy-Pad\'{e} approximations.

The rest of this paper is arranged as follows:  Section \ref{sec2}  introduces a novel the fractional order of rational Bessel functions (FRB). Section \ref{sec3} describes a brief formulation of quasilinearization method (QLM) introduced by \cite{refff05}. In section \ref{sec4} at first, by utilizing QLM over Thomas-Fermi equation  a sequence of linear differential equations are obtained then in each iteration solve the linear differential equation by the fractional order of rational Bessel functions collocation method (FRBC) that we name FRBC-QLM for used method. We  in  section \ref{sec5} compare our solutions with some well-known results, comparisons show that the present solutions are highly accurate, also we describe our results via tables and figures. Finally we give a brief conclusion in section \ref{sec6}.

\section{The fractional order of rational Bessel functions $(FRB)$}\label{sec2}
The Bessel functions arise in many problems in physics possessing cylindrical symmetry, such as the vibrations of circular drumheads and the radial modes in optical fibers. Bessel functions are usually defined as a particular solution of a linear differential equation of the second order which known as Bessel's equation. Bessel functions first defined by the Daniel Bernoulli on heavy chains (1738) and then generalized by Friedrich Bessel. More general Bessel functions were studied by Leonhard Euler in (1781) and in his study of the vibrating membrane in (1764) \cite{reffff01,reffff02}.
\subsection{Definition of Bessel polynomials}

Bessel differential equation of order $n\in\mathbb{R}$ is:
\begin{eqnarray}\label{eqqq04}
x^2\frac{d^2y(x)}{dx^2}+x\frac{dy(x)}{dx}+(x^2-n^2) y(x)=0,~~x\in(-\infty,\infty).
\end{eqnarray}

One of the solutions of equation (\ref{eqqq04}) by applying the method of Frobenius as follows \cite{reffff03}:
\begin{eqnarray}\label{eqqq05}
J_{n}(x)=\sum_{r=0}^{ \infty}\frac{(-1)^r}{r!(n+r+1)!}(\frac{x}{2})^{2r+n},
\end{eqnarray}
where series (\ref{eqqq05})  is convergent for all $x\in(-\infty,\infty)$.

Bessel functions and polynomials are used to solve the more number of problems in physics, engineering, mathematics, and etc., such as Blasius equation, Lane-Emden equations, integro-differential equations of the fractional order, unsteady gas equation, systems of linear Volterra integral equations, high-order linear complex differential equations in circular domains, systems of high-order linear Fredholm integro-differential equations, etc. \cite{reffff04,reffff05,reffff06,reffff07,reffff08,reffff09,reffff10,reffff11,reffff12,reffff13,reffff1301,reffff13011,reffff1302}.

Bessel polynomials has been introduced as follows \cite{reffff12}:
\begin{eqnarray}\label{eqqq06}
B_{n}(x)=\sum_{r=0}^{[\frac{N-n}{2}]}\frac{(-1)^r}{r!(n+r)!}(\frac{x}{2})^{2r+n},~~x\in[0,1].
\end{eqnarray}
where $n\in\mathbb{N}$, and $N$ is the number of basis of Bessel polynomials.

Let us define $\Gamma=\{x|~0\leq x \leq 1 \}$ and \\
$L^{2}_{w}(\Gamma)=\{~v :\Gamma \rightarrow \mathbb{R}| v$ is measurable and $\parallel v \parallel_{w} < \infty \}$, where
$$\parallel v \parallel_{w}=\left(\int^{1}_{0}|v(x)|^{2}w(x)dx\right)^{1/2},$$
with $w(x)=1$, is the norm induced by inner product of the space $L^{2}_{w}(\Gamma)$ as follows:
$$\langle v,u\rangle_{w}=\int^{1}_{0}{v(x)u(x)w(x)}dx.$$

Now, suppose that\\
$$\mathfrak{B}=~span\{B_{0}, B_{1},\dots, B_{N}\},$$
where $\mathfrak{B}$ is a finite-dimensional subspace of $L^{2}_{w}(\Gamma)$, dim $\mathfrak{B} = N+1$, so $\mathfrak{B}$ is a closed subspace of $L^{2}(\Gamma)$. Therefore, $\mathfrak{B}$ is a complete subspace of $L^{2}(\Gamma)$. Assume that f be an arbitrary element in $L^{2}(\Gamma)$. Thus f has a unique best approximation in $\mathfrak{B}$ subspace, say $\hat{b}\in \mathfrak{B}$, that is
\begin{eqnarray}
\nonumber\exists~ \hat{b}\in\mathfrak{B}, ~~~ \forall ~b\in \mathfrak{B},~~\parallel f-\hat{b}\parallel \leq \parallel f-b\parallel.
\end{eqnarray}
Notice that we can write $b$ vector as a combination of the basis vectors of $\mathfrak{B}$ subspace. 

We know function of $f$  can be expanded by $N+1$ terms of Bessel polynomials as: 
\begin{eqnarray}
f=f_{N}+R,
\nonumber\end{eqnarray}
that is
\begin{equation}\label{eqqq07}
f_{N}=\sum^{N}_{n=0}{a_{n}B_{n}(x)}=A^{T}B(x),
\end{equation}
where $B(x)=[B_{0}(x), B_{1}(x),\cdots, B_{N}(x)]^{T}$ and $R\in\mathfrak{B}^{\perp}$ that ${\mathfrak B}^{\perp}$ is the orthogonal complement. So  $f-f_{N}\in\mathfrak{B}^{\perp}$ and $b\in\mathfrak{B}$ are orthogonal which we denote it by:
\begin{eqnarray}
\nonumber f-f_{N}\perp b,
\end{eqnarray}
thus $f-f_{N}$ vector is orthogonal over all of basis vectors of $\mathfrak{B}$ subspace as:
\begin{eqnarray}
\nonumber \langle f-f_{N},B_{i}(x)\rangle_{w}=\langle f-A^{T}B(x),B_{i}(x)\rangle_{w}=0,~i=0, 1,\cdots, N,
\end{eqnarray}
hence
\begin{eqnarray}
&&\nonumber\langle f-A^{T}B(x),B^{T}(x)\rangle_{w}=0,
\end{eqnarray}
therefore A can be obtained by
\begin{eqnarray}
&&\nonumber\langle f,B^{T}(x)\rangle_{w}=\langle A^{T}B(x),B^{T}(x)\rangle_{w},\\
\nonumber\\
\nonumber &&A^{T}=\langle f,B^{T}(x)\rangle_{w}\langle B(x),B^{T}(x)\rangle_{w}^{-1},~n=0, 1,\cdots, N.
\end{eqnarray}

\subsection{Definition of FRB}
Some researchers have proposed the series expansions  $\sum_{i=0}^{N}{c_{i}x^{i\alpha}},~(\alpha>0)$ to solve the fractional differential equations, for instance, Bhrawy et al have constructed shifted fractional-order Jacobi orthogonal functions  to solve the nonlinear initial value problem of fractional order $\alpha$ and a class of time-fractional partial differential equations with variable coefficient \cite{reffff14,reffff15}. Authors \cite{reffff16,reffff17} have proposed  fractional-order Legendre functions to solve fractional-order differential equations and the time-fractional convection-diffusion equation. Alshbool et al. have utilized operational matrices of new fractional Bernstein functions for approximating solutions to fractional differential equations \cite{reffff18}. Parand and Delkhosh have introduced the fractional order of the Chebyshev functions for solving Volterra's population growth model of arbitrary order \cite{reffff188}.

We have applied new FRB to solve the Thomas-Fermi equation in semi-infinite interval, because by analytical solution for this problem has proved that the answer is as fractional forms \cite{refd07}. Now we define FRB,  $\{ FB_{n}\}$ by:
\begin{eqnarray}
FB_{n}^{\alpha}(x,L)=B_{n}(\frac {x^{\alpha}}{x^{\alpha} + L})~,~~n=0, 1,\cdots,N \nonumber
\end{eqnarray}
or
\begin{eqnarray}\label{eqqq09}
FB_{n}^{\alpha}(x,L)=\sum_{r=0}^{[\frac{N-n}{2}]}\frac{(-1)^r}{r!(n+r)!}(\frac{x^{\alpha}}{2(x^{\alpha} + L)})^{2r+n},~~n=0, 1,\cdots,N
\end{eqnarray}
where $\alpha>0$, $x\in [0,\infty)$, $B_{n}(x)$ is Bessel polynomials of order $n$, and the constant parameter $L>0$ is a scaling/stretching factor.

Let us define $\Lambda =\{x|~0\leq x < \infty \}$ and \\
$L^{2}_{w}(\Lambda )=\{~z :\Lambda \rightarrow \mathbb{R}| z$ is measurable and $\parallel z \parallel_{w} < \infty \}$, where
$$\parallel z \parallel_{w}=\left(\int^{\infty}_{0}|z(x)|^{2}w(x,L)dx\right)^{1/2},$$
with $w(x,L)=\frac{\alpha x^{\alpha-1}L}{(x^{\alpha}+L)^{2}}$, is the norm induced by inner product of the space $L^{2}_{w}(\Lambda)$ as follows:
$$\langle z,g\rangle_{w}=\int^{\infty}_{0}{z(x)g(x)w(x,L)}dx.$$

Now, suppose that\\
$$\mathfrak{FB}=~span\{FB^{\alpha}_{0}(x,L), FB^{\alpha}_{1}(x,L),\dots, FB^{\alpha}_{N}(x,L)\},$$
Let $y\in L^{2}(\Gamma)$ be a function defined over interval $[0,\infty)$ can be expanded by $N+1$ terms of FRB as: 
\begin{eqnarray}\label{eqqq10}
y_{N}(x)=\sum^{N}_{n=0}{a_{n}FB^{\alpha}_{n}(x,L)}=A^{T}FB(x,L),
\end{eqnarray}
where $FB(x,L)=[FB^{\alpha}_{0}(x,L), FB^{\alpha}_{1}(x,L),\dots, FB^{\alpha}_{N}(x,L)]$. Hence
\begin{eqnarray}\label{eqqq11}
&&\nonumber\langle y-A^{T} FB(x,L), FB^{T}(x,L)\rangle_{w}=0,
\end{eqnarray}
therefore A can be obtained by
\begin{eqnarray}
&&\nonumber\langle y, FB^{T}(x,L)\rangle_{w}=\langle A^{T} FB(x,L), FB^{T}(x,L)\rangle_{w},\\
\nonumber\\
\nonumber&&A^{T}=\langle y, FB^{T}(x,L)\rangle_{w}\langle  FB(x,L), FB^{T}(x,L)\rangle_{w}^{-1},~n=0, 1,\cdots, N.
\end{eqnarray}

\section{The quasilinearization method (QLM)}\label{sec3}
The QLM is a generalization of the Newton-Raphson method \cite{refff01,refff02} to solve the nonlinear differential equation as a limit of approximating the nonlinear terms by an iterative sequence of linear expressions. Bellman and Kalaba was originally introduced the QLM method about fifty years ago \cite{refff03,refff04} . The QLM techniques are based on the linearization of the high order ordinary/partial differential equation and require the solution of a linear ordinary differential equation at each iteration. Mandelzweig and Tabakin \cite{refff05} have determined general conditions for the quadratic, monotonic and uniform convergence of the QLM method to solve both initial and boundary value problems in nonlinear ordinary $n$th order differential equation in $N$-dimensional space. Recently, QLM method has been successfully applied by researchers to solve various type of fractional differential equations and some ordinary nonlinear equation \cite{refff06,refff07,refff08,refff09}.

We have been supposed second-order nonlinear ordinary differential equation in one variable on the interval $[0, \infty)$ as follows:
\begin{eqnarray}\label{eqqq12}
\frac{d^2u}{dx^2}=F(u'(x),u(x),x),
\end{eqnarray}
with the boundary conditions: $u(0)=A,~u(\infty)=B$, where $A$ and $B$ are real constants and $F$ is nonlinear functions.

By using the QLM for solving Eq. (\ref{eqqq11}) determines the $(r+1)$th iterative approximation $u_{r+1}(t)$ as a solution of the linear differential equation:
\begin{equation}\label{eqqq12}
\frac{d^2u_{r+1}}{dx^2}=F(u'_r,u_r,x)+(u_{r+1}-u_r)F_u(u'_r,u_r,x)+(u'_{r+1}-u'_r)F_{u'}(u'_r,u_r,x),
\end{equation}
with the boundary conditions:
\begin{equation}\label{eqqq13}
u_{r+1}(0)=A,~~~~~u_{r+1}(\infty)=B,
\end{equation}
where $~r=0,1,2,\cdots$ and the functions $F_u = \partial F/ \partial u$ and $F_{u'} = \partial F / \partial u'$ are functional derivatives of functional $F(u'_r,u_r,x)$.

\section{Solution of Thomas-Fermi equation by FRBC-QLM}\label{sec4}
By utilizing QLM technique on Eq. (\ref{eqqq01}), we have
\begin{eqnarray}\label{eqqq14}
\frac{d^{2}y_{r+1}(x)}{dx^2} - \frac{3}{2\sqrt{x}}(y_{r}(x))^{1/2}y_{r+1}(x)=-\frac{1}{2\sqrt{x}}(y_{r}(x))^{3/2},
\end{eqnarray}
with the boundary conditions:
\begin{eqnarray}\label{eqqq15}
y_{r+1}(0)=1,~~~~~~y_{r+1}(\infty)=0,
\end{eqnarray}
where $r = 0, 1, 2,\cdots .$\\

For rapid convergence is actually enough that the initial guess is sufficiently good to ensure the smallness of just one of the quantity $q_{r} = k||y_{r+1} - y_{r}||$, where k is a constant independent of r. Usually, it is advantageous that $y_{0}(t)$
would satisfy at least one of the boundary conditions Eq.  (\ref{eqqq15}) \cite{refff06}, thus set
 $y_{0}(x)=1$ for the initial guess of Thomas-Fermi equation. In this paper have been considered two terms $\frac{1}{x^2+1}$ and $\frac{x}{x^2+1}$ to satisfy boundary conditions Eq. (\ref{eqqq15}). Thus we can approximate $y_{r+1}(x)$  by $N+1$ basis of FRB as:

\begin{equation}\label{eqqq16}
y_{r+1}(x) \thickapprox y_{N,r+1}(x)=\frac{1}{x^2+1}+\frac{x}{x^2+1}\sum^{N}_{n=0}{\hat{c_{i}}FB^{\alpha}_{n}(x,L)}.
\end{equation}
where $\alpha>0$ and $r = 0, 1, 2,\cdots .$ In all of the spectral methods, the purpose is to find $\hat{c_{i}}$ coefficients.\\
To apply the collocation method, we have constructed the residual function for  $(r+1)$th iteration in QLM method by substituting $y_{r+1}(x)$ by $y_{N,r+1}(x)$ into Eq. (\ref{eqqq14}) as following:
\begin{equation}\label{eqqq17}
Res_{r+1}(x)=\frac{d^{2}y_{N,r+1}}{dx^{2}}-\frac{3}{2\sqrt{x}}(y_{r}(x))^{1/2}y_{N,r+1}(x)+\frac{1}{2\sqrt{x}}(y_{r}(x))^{3/2}.
\end{equation}
 A method for forcing the residual function (\ref{eqqq17}) to zero can be defined as collocation algorithm. There is no limitation to choose the point in collocation method. The $N+1$ collocation points which are roots of rational Chebyshev functions on interval [0,$\infty$) (i.e.$~~x_{i}=(1-cos(\frac{(2i-1)\Pi}{2N+2}))/(1+cos(\frac{(2i-1)\Pi}{2N+2})),~i=1, 2,\dots, N+1$ \cite{reff05}) have been substituted $Res_{r+1}(x)$, therefore:
\begin{equation}\label{eqqq18}
Res_{r+1}(x_{i})=0,~~i=0, 1, ,\cdots, N+1.
\end{equation}
A nonlinear system of equations has been obtained, all of these equations can be solved by Newton method for the unknown coefficients. We have also done all of computations by Maple 2015 on PC with CPU Core i5, Windows 7 64bit, and 8GB of RAM.

Now we can employ the FRBC-QLM  iterative algorithm to solve Thomas-Fermi equation as follows:\\\\
\begin{bf}
%\large
BEGIN
\end{bf}
\begin{enumerate}
\item[\begin{bf}(1)\end{bf}] Input variable of $I$ that is the number of iterations of QLM method.
\item[\begin{bf}(2)\end{bf}] Input variable of $N$ that is the number of basic of the FRB.
\item[\begin{bf}(3)\end{bf}]  Set $y_{N,0}(x)=1$.
\item[\begin{bf} For \end{bf}]r=0 \begin{bf} to \end{bf} $I$ \textbf{do} 
\begin{enumerate}
\item[\begin{bf}(4)\end{bf}]  Construct the series (\ref{eqqq16}) for approximating $y_{r+1}(x)$ as $y_{N,r+1}(x)$.
\item[\begin{bf}(5)\end{bf}] Construct  the  linear differential equation (\ref{eqqq17}) by using QLM method on the Eq. (\ref{eqqq01}).
\item[\begin{bf}(6)\end{bf}] Substitute $y_{N,r+1}(x)$ into the  equation (\ref{eqqq17}) and create residual function $Res_{r+1}(x)$.\\ Now we have $N +1$ unknown $\{ \hat{c_{i}}\}_{0}^{N}$. To obtain these unknown coefficients, we need $N +1$ equations.
\item[\begin{bf}(7)\end{bf}] Choose the roots of order $N+1$ of Rational Chebyshev functions as $N+1$ collocation points: $\{ x_{i}\}_{0}^{N}$.
\item[\begin{bf}(8)\end{bf}] Substitute collocation points $\{ x_{i}\}_{0}^{N}$ into the $Res_{r+1}(x)$ and create the $N+1$ equations.
\item[\begin{bf}(9)\end{bf}] Solve the $N+1$ linear equations with $N+1$ unknown coefficients, for calculate $y_{N,r+1}(x)$.
\end{enumerate}
\item[\begin{bf}End\end{bf}]  \begin{bf} For \end{bf}
\end{enumerate}
\begin{bf}
%\large
END
\end{bf}

\section{Numerical Results}\label{sec5}
The initial slope $y'(0)$ is difficult to compute by any means, and plays an important role in determining many physical properties of the Thomas-Fermi atom. It determines the energy of a neutral atom in the Thomas-Fermi approximation. Zaitsev et al \cite{refh6} have showed that methods of Runge-Kutta and Adams-Bashforth can apply to solve the Thomas-Fermi equation in semi-infinite interval, although their methods are ill-condition and have not high accuracy for more scheme. Exact solution for Thomas-Fermi differential equation, which is defined in the semi-infinite interval and has singularity at $x=0$ and its boundary condition occurs at infinity, is not available, so approximating this solution is very important.

Table 1 shown a list of the number of calculations y'(0) of the Thomas-Fermi potential. As can be seen, some researchers have achieved good results and accurate. The last three rows show best approximations of $y'(0)$ for various value of $N$ and fixed value of $L=1$ by present method which shows that the present solution is highly accurate. Tables 2 and 3 show values obtained of $y(x)$ and $y'(x)$ by present method respectively, for different values of $N$ and iteration=45. Obviously, Table 4 and 5 have presented some numerical example to illustrate the accuracy and convergence of our suggested method by increasing the number of points and iterations. It should be mentioned that all calculations are done by software Maple for various values $N$ and iterations. Fig. 1  shows the resulting graph of Thomas-Fermi equation  graph obtained by present method for $N=200$ and iteration 45 which tends to zero as $x$ increases by boundary condition $y(\infty) = 0$, and graphs of residual error of the problem with $N=50,100,150,200$, and iteration 45, note that the residual error decreases with the increase of the collocation points. Comparing the computed results by this method with the others shows that this method provides more accurate and numerically stable solutions than those obtained by other methods.

\section{Conclusion}\label{sec6}
The fundamental goal of this paper has been to construct an approximation to the solution of nonlinear Thomas-Fermi equation in a semi-infinite interval which has singularity at $x = 0$ and its boundary condition occurred in infinity. In the above discussion, we applied a new method to solve the Thomas-Fermi equation that is nonlinear ordinary differential equation on a semi-infinite interval. By using an analytical method for solving Thomas-Fermi equation has proved that the answer of this problem is as fractional forms \cite{refd07}. So for the first time, we solved the problem based on the new fractional order of rational Bessel functions  without any domain truncation or transformation of the domain of the problem to a finite domain. In this work first, by utilizing QLM over Thomas-Fermi equation a sequence of linear differential equations are obtained. Second, in each iteration solve the linear differential equation by a novel FRBC method. We obtained accurate to 30 decimal places for initial slope, $y'(0) = -1.588071022611375312718684509423$, only by using 200 collocation points and successfully have been applied to find the most accurate values of $y(x)$ and $y'(x)$. A known open problem in spectral methods is finding the optimal value for $L$ \cite{reff001}, but in this paper, we do not need to calculate the optimal value for $L$ so that, we set  $L = 1$. The numerical results of solving this problem show that this method is  higher accurate than obtained results of other famous methods. Finally, comparison results have shown that the present method is an acceptable approach and good candidate to solve this type of problems that occur in semi-infinite interval and the nonlinear singular two point boundary value problems effectively.

%%%%%%

%%%%%%%%%%%%%%%%%%%%%%%%%%%%%%%%%%%%%%%%%%

\begin{figure}[!ht]
\centering
\subfigure[Graphs of residual error.]{
\includegraphics*[width=5.5cm]{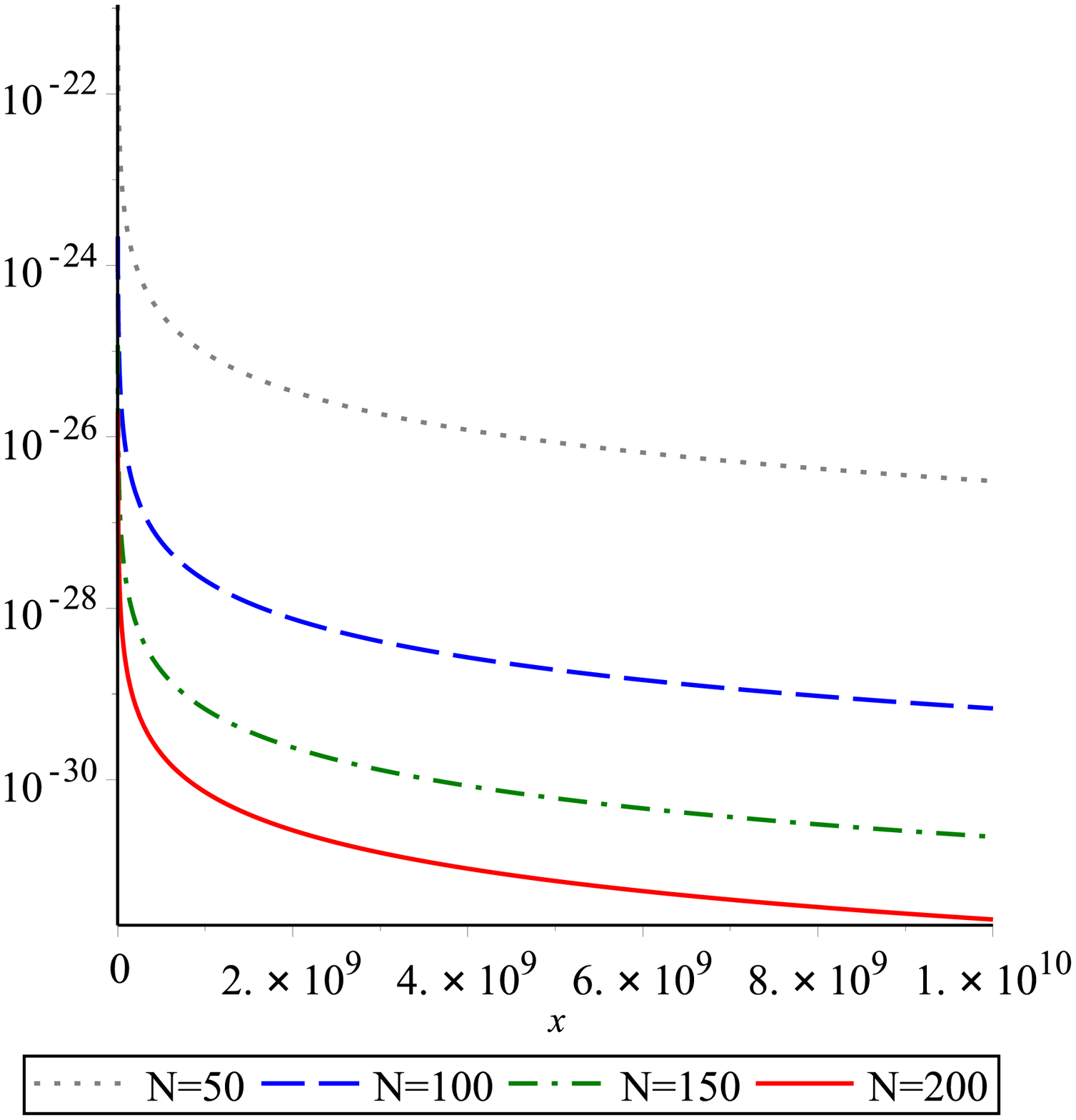}}
\hspace{1mm}
\subfigure[Graph of $y(x)$.]{
\includegraphics*[width=5.5cm]{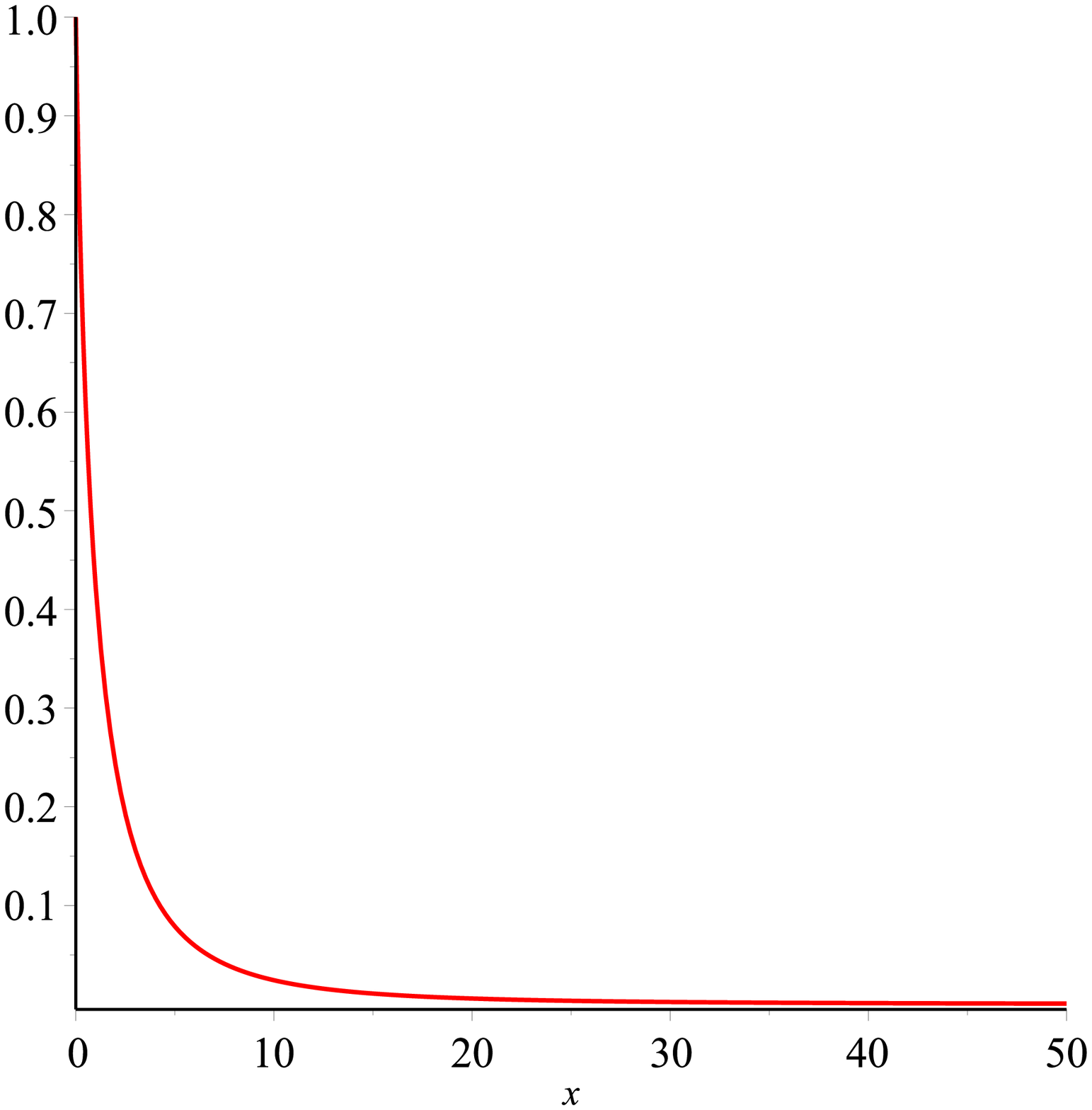}}
\caption{Graphs of residual error with $N=50,100,150,200$, and iteration 45, and Thomas-Fermi graph obtained by present method.}
\end{figure}

%%%%%%%%%%%%%%%%%%%%%%%%%%%%%%%%%%%%%%%%%%

%%%%%%%%%%%%%%%%%%%%%%%%%%%%%%%%%%%%%%%%%%
\footnotesize 
\begin{center}
\begin{table}[ht!]
\caption{Comparison of obtained values of $y'(0)$ by researchers, inaccurate digits are in bold
face.}
\centering
%\footnotesize 
\begin{tabular}[ht!]{ll}
\toprule
  Author/Authors & Obtained value of $y'(0)$  \\ 
\midrule
   Fermi (1928) \cite{refd05} & -1.58     \\
   Baker (1930) \cite{refd07}  & -1.588\textbf{558}     \\
   Bush and Caldwell (1931) \cite{reff11} & -1.58\textbf{9}     \\
   Miranda (1934) \cite{reff12}  & -1.5880\textbf{464}     \\
   Slater and Krutter (1935) \cite{refd03}  & -1.5880\textbf{8}     \\
   Feynman et al (1949) \cite{refd04} & -1.588\textbf{75}     \\
   Kobayashi et al. (1955) \cite{refd31}  & -1.58807\textbf{0972}     \\
   Mason (1964) \cite{refd09} & -1.5880710     \\
   Laurenzi (1990) \cite{refd06} & -1.588\textbf{588}     \\
   MacLeod (1992) \cite{reff03} & -1.5880710226     \\
   Wazwaz (1999) \cite{refh3} & -1.58807\textbf{6779}     \\
   Epele et al (1999) \cite{refm18} & -1.588\textbf{1}     \\
   Esposito (2002) \cite{refd14} & -1.588     \\
   Liao (2003) \cite{refd16} & -1.58\textbf{712}     \\
   Khan and Xu (2007) \cite{refm19} & -1.58\textbf{6494973}     \\
   El-Nahhas (2008) \cite{refd18} & -1.5\textbf{5167}     \\
   Yao (2008) \cite{refm12} & -1.5880\textbf{04950}     \\
   Fernandez (2008) \cite{refm17} & -1.588071022611375313     \\
   Parand and Shahini (2009) \cite{reff05} & -1.58807\textbf{02966}     \\
   Marinca and Herianu (2011) \cite{refh4} & -1.5880\textbf{659888}     \\
   Oulne (2011) \cite{refh7} & -1.5880710\textbf{34}     \\
   Abbasbandy and Bervillier (2011) \cite{reff10} & -1.588071022611375312718\textbf{9} \\
   Zhu et al. (2012) \cite{refh5} & -1.58\textbf{794357}     \\
   Turkylmazoglu (2012) \cite{refm11} & -1.5880\textbf{1}     \\
   Zhao et al (2012) \cite{refm13} & -1.5880710226     \\
   Parand et al (2013) \cite{reff09} & -1.58807\textbf{0339}     \\
   Boyd (2013) (with m=600) \cite{reff01} & -1.5880710226113753127186845     \\
   Tavassoli Kajan et al. (2013) \cite{reff04} & -1.58807102261137\textbf{4}     \\
   Amore et al (2014) \cite{refm16} & -1.588071022611375312718684508     \\
   Marinca and Ene (2014) \cite{refd25} & -1.588071\textbf{9992}     \\
   Bayatbabolghani and Parand(2014)\cite{reff07} & -1.588071     \\
   Kilicman et al (2014) \cite{reff02} & -1.588071\textbf{347}     \\
   Liu and Zhu (2015) \cite{reff08} & -1.58807\textbf{2}     \\
  Parand et al [this article, N=100]&-1.5880710226113753127   \\
'' ''[N=150]& -1.5880710226113753127186845\\
'' ''[N=200]& -1.588071022611375312718684509423  \\
\bottomrule
\end{tabular}
\end{table}
\end{center}
%\normalsize
%%%%%%%%%%%%%%%%%%%%%%%%%%%%%%%%%%%%%%%%%

%%%%%%%%%%%%%%%%%%%%%%%%%%%%%%%%%%%%%%%%%%
\footnotesize 
\begin{center}
\begin{table}[ht!]
\caption{Values of $y(x)$ for various values of $x$ with iteration 45 and N=200 }
\centering
%\footnotesize 
\begin{tabular}[ht!]{lc||lc}
\toprule
   $x$  & $y(x)$ & $x$ & $y(x)$    \\ 
\midrule       
   0.25   & 0.755201465313331276073659062048  & 30     &  0.002255836616202855884224076079       \\ 
    0.50  & 0.606986383355979909494446070174  & 40     &  0.001113635638833368812571167838   \\ 
    0.75  & 0.502346846412368627446521794036  & 50     &  0.000632254782984904726779778752     \\ 
    1.00  & 0.424008052080705600224612007418  & 60     &  0.000393911366685417002026602841    \\ 
    1.25  & 0.363201414459514114681451617277  & 70     &  0.000262265299812011293741733594       \\ 
    1.50  & 0.314777463700458172973580939810  & 80     &  0.000183545759740710240155439613          \\ 
    1.75  & 0.275451327996091785680113852782  & 90     &  0.000133545828953734623590590348       \\ 
    2.00  & 0.243008507161119555299806749733  & 100   &  0.000100242568139407331685585932          \\ 
    2.25  & 0.215894626576130144431137496637  & 200   &  0.000014501803496945764680403612          \\ 
    2.50  & 0.192984123458000701287136925252  & 300   &  0.000004548571953616680184257373            \\ 
    2.75  & 0.173441292490063451179594691770  & 400   &  0.000001979732628112504742575796\\ 
    3.00  & 0.156632673216495841339813440477  & 500   &  0.000001034077168199939706035333           \\ 
    3.25  & 0.142069642692650781317847467819  & 600   &  0.000000606868769667525133710573    \\
    3.5    & 0.129369596993799111381550504704  & 700   &  0.000000386176515703798616221824     \\
    3.75  & 0.118229001616846686506107664622  & 800   &  0.000000260813730499833635331580     \\
    4.00  & 0.108404256918907711089847680321  & 900   &  0.000000184372415135065178917921     \\
    4.25  & 0.099697845864740046922595440377  & 1000 &  0.000000135127477354105831539440   \\
    4.50  & 0.091948133826563845114632113892  & 2000 &  0.000000017339847516138509108205   \\
    4.75  & 0.085021743728059499429131697623  & 3000 &  0.000000005189408334543857341875     \\
    5.00  & 0.078807779251369904256091892542  & 4000 &  0.000000002201209082423027362721    \\
    6.00  & 0.059422949250422580797949567059  & 5000 &  0.000000001130926706419984771574     \\
    7.00  & 0.046097818604498589876456260102  & 6000 &   0.000000000656056637887703224451    \\
    8.00  & 0.036587255264676802392315804375  & 7000 &  0.000000000413886522087042381058    \\
    9.00  & 0.029590935270546873724362041806  & 8000 &  0.000000000277658195353364145555     \\
  10.00  & 0.024314292988680864190110388176  & 9000 &  0.000000000195225879692802747067   \\
   20.00 & 0.005784941191566940442010571504  &10000&  0.000000000142450044462688615523   \\

\bottomrule
\end{tabular}
\end{table}
\end{center}
%%%%%%%%%%%%%%%%%%%%%%%%%%%%%%%%%%%%%%%%%%

%%%%%%%%%%%%%%%%%%%%%%%%%%%%%%%%%%%%%%%%%%
\footnotesize 
\begin{center}
\begin{table}[ht!]
\caption{Values of $y'(x)$ for various values of $x$ with iteration 45 and $N=200$}
\centering
%\footnotesize 
\begin{tabular}[ht!]{lc||lc}
\toprule
   $x$  & $y'(x)$ & $x$ & $y'(x)$    \\ 
\midrule
    0.25  &- 0.722306984910234919519668083864  & 30     &  -0.000180670006476992635016139693       \\ 
    0.50  & -0.489411612574538088647005847557  & 40     &  -0.000069668028540325866316070590   \\ 
    0.75  & -0.358306880167513621987250767311  & 50     &  -0.000032498902048258814624200679     \\ 
    1.00  & -0.273989051593306251989464686519  & 60     &  -0.000017197700083099862592575862    \\ 
    1.25  & -0.215794130300733601274300529492  & 70     &  -0.000009956533393052361268495181       \\ 
    1.50  & -0.173738799013945185681936465228  & 80     &  -0.000006166195528764075475304167          \\ 
    1.75  & -0.142320937196893658885960452965  & 90     &  -0.000004024473703766734693502284     \\ 
    2.00  & -0.118243191625487620571255867534  & 100   &  -0.000002739351068678330086744132          \\ 
    2.25  & -0.099409321201447030009355248089  & 200   &  -0.000000205753231647526892605535          \\ 
    2.50  & -0.084426186798809043812545918214  & 300   &  -0.000000043659496185302904545958          \\ 
    2.75  & -0.072335044846097235621822505698  & 400   &  -0.000000014366823059961810213641     \\ 
    3.00  & -0.062457130854120976228704899999  & 500   &  -0.000000006034363442475256759610            \\ 
    3.25  & -0.054300422911798016711579461695  & 600   &  -0.000000002961822515102276227611    \\
    3.5    & -0.047501046582295208950097689053  & 700   &  -0.000000001619832187577198029687     \\
    3.75  & -0.041785207716826396797995883443  & 800   &  -0.000000000959243855994648160301     \\
    4.00  & -0.036943757824123486354813738987  & 900   &  -0.000000000603766177055436240659    \\
    4.25  & -0.032814785443993540385309640573  & 1000 &  -0.000000000398801070822799359788 \\
    4.50  & -0.029271448448803843379269452384  & 2000 &  -0.000000000025760853699197105401   \\
    4.75  & -0.026213311168397937715703460851  & 3000 &  -0.000000000005153001176402320032     \\
    5.00  & -0.023560074954700512881180449080  & 4000 &  -0.000000000001641618607042605687    \\
    6.00  & -0.015867549533407079812737615662  & 5000 &  -0.000000000000675339712187831946      \\
    7.00  & -0.011142531814867088405578460800  & 6000 &   -0.000000000000326676807336313998     \\
    8.00  & -0.008088602969645474322126751847  & 7000 &   -0.000000000000176730816571737264    \\
    9.00  & -0.006033074714457392439143608348  & 8000 &   -0.000000000000103777992089730866    \\
  10.00  & -0.004602881871269254502543511851  & 9000 &  -0.0000000000000648790244833915206   \\
   20.00 & -0.000647254332777692033047085589  &10000&  -0.0000000000000426161649603483992   \\

\bottomrule
\end{tabular}
\end{table}
\end{center}
%%%%%%%%%%%%%%%%%%%%%%%%%%%%%%%%%%%%%%%%%%

%%%%%%%%%%%%%%%%%%%%%%%%%%%%%%%%%%%%%%%
\begin{center}
\begin{table}[h!]
\center
\caption{Numerical results of solution $y(x)$ with various values of $x$, $N$ and iterations.  }
\begin{tabular}{| m{1em} | m{1.2em} | m{11em} m{16em} m{16em} | } 
\hline
$~n$ & $x$ &~~~~~~~~~15th iteration & ~~~~~~~~~~~~~~~30th iteration & ~~~~~~~~~~~~~~~45th iteration\\
\hline % n=50
50 &10 100 200 300 400 500 & 
0.02431429298880024169       			    0.00010024257693932961			       0.00001450155426224488 				  0.00000454855611182438			     0.00000198237212867001			 0.00000104219493450734 & 
0.024314292988680868310550970861		    0.000100242568139360181452116145	       0.000014501803498782770002228872   0.000004548571957198050067791132   0.000001979732627261138334433993	0.000001034077157264396322609605 & 
0.024314292988680865622793862360   	    0.000100242568139361977721849610   	       0.000014501803498894699818187404   0.000004548571957423423877537837   0.000001979732627197689616236145   	0.000001034077156421974180242224\\
\hline
\hline % n=100
100 & 10 100 200 300 400 500 & 
0.02431429298868156639    0.00010024257310248919   0.00001450191210978454    0.00000454926463044528     0.00000198229714272620 			  0.00000104078259214971& 
0.024314292988680864190110392994		   0.000100242568139407331732468387      0.000014501803496945767789317745       0.000004548571953616664783487078     0.000001979732628112388963029981		   0.000001034077168200029689420358 & 
0.024314292988680864190110392609 	   0.000100242568139407331736524506      0.000014501803496945768054100507       0.000004548571953616663090340995     0.000001979732628112377070251772  		   0.000001034077168200025910092884 \\ 
\hline
\hline % n=150
150& 10 100 200 300 400 500 & 
0.02431429298868219619    0.00010024257760384483    0.00001450201192996564       0.00000454989194540830      0.00000198461606782965      0.00000104695705516519  & 
0.024314292988680864190110388158	           0.000100242568139407331685519490      0.000014501803496945764679973963     0.000004548571953616680174309015        0.000001979732628112504877013189    0.000001034077168199940307972145& 
0.024314292988680864190110388161            0.000100242568139407331685518495      0.000014501803496945764680397208 	0.000004548571953616680172305663       0.000001979732628112504797936584     0.000001034077168199940285585132\\
\hline
\hline %n=200
200&10 100 200 300 400 500 & 
0.02431429298868235435    0.00010024257868296334     0.00001450203484084579   0.00000455003833574595   0.00000198515363835884  0.00000104836091033391  & 
0.024314292988680864190110388176   0.000100242568139407331685586056   0.000014501803496945764680404365  0.000004548571953616680184259791    0.000001979732628112504742514001    0.000001034077168199939706184215& 
0.024314292988680864190110388176      0.000100242568139407331685585932   0.000014501803496945764680403612    0.000004548571953616680184257373     0.000001979732628112504742575794   0.000001034077168199939706035334  \\ 
\hline
\end{tabular}
\end{table}
\end{center}
%%%%%%%%%%%%%%%%%%%%%%%%%%%%%%%%%%%%%%%%%%

%%%%%%%%%%%%%%%%%%%%%%%%%%%%%%%%%%%%%%%
\begin{center}
\begin{table}[h!]
\center
\caption{Numerical results of solution $y'(x)$ with various values of $x$, $N$ and iterations.  }
\begin{tabular}{| m{1em} | m{1.2em} | m{11.3em} m{16.1em} m{16.1em} | } 
\hline
$~n$ & $x$ &~~~~~~~~~15th iteration & ~~~~~~~~~~~~~~~30th iteration & ~~~~~~~~~~~~~~~45th iteration\\
\hline % n=50
50 &0 10 100 200 300 400 500 & 
-1.58798412034597225837   -0.00460288187129760549	 -0.00000273935089878581  -0.00000020575729375812	      -0.00000004364704984167	  	-0.00000001432569675626 	-0.00000000596694589883&  
-1.588071024773768750505101814302    -0.004602881871269255485186656747  -0.000002739351068678941009467810    -0.000000205753231614201535153310 	-0.000000043659496195020703717486  -0.000000014366823136534945762465 -0.000000006034363563360980464672 & 
-1.588071022461220318498896590154  -0.004602881871269254843810395419   -0.000002739351068678834237868374  -0.000000205753231612157754759952  -0.000000043659496195672746505803  -0.000000014366823142132220490812 -0.000000006034363572464639504768\\
\hline
\hline % n=100
100 &0 10 100 200 300 400 500 & 
-1.58806849943926523283   -0.00460288187126903604   -0.00000273935085186272   -0.00000020575078049313  -0.00000004364886726607   -0.00000001433829371433   -0.00000000597884939751& 
-1.588071022611375312713595678174   -0.004602881871269254502543510274   -0.000002739351068678330084199975   -0.000000205753231647526767499279   -0.000000043659496185303706721008   -0.000000014366823059962243265195 -0.000000006034363442469847608554 & 
-1.588071022611375312724621425500  -0.004602881871269254502543510545   -0.000002739351068678330084039598      -0.000000205753231647526755494256    -0.000000043659496185303784906289     -0.000000014366823059962315807005 -0.000000006034363442469561014542\\ 
\hline
\hline % n=150
150&0 10 100 200 300 400 500 & 
-1.58807102261138081141    -0.00460288187126883975  -0.00000273935065143065   -0.00000020574848020790   -0.00000004365807090310 -0.00000001431185061144  -0.00000000592789998809 & 
-1.588071022611375312718684511366   -0.004602881871269254502543511857  -0.000002739351068678330086740561  -0.000000205753231647526892507632  -0.000000043659496185302905154071  -0.000000014366823059961806789804  -0.000000006034363442475251992346 & 
-1.588071022611375312718684517975  -0.004602881871269254502543511856   -0.000002739351068678330086729788 	 -0.000000205753231647526892514803   -0.000000043659496185302905165385  -0.000000014366823059961806903847 -0.000000006034363442475252124155 \\
\hline
\hline %n=200
200 &0 10 100 200 300 400 500 & 
-1.58807102261137532732     -0.00460288187126879049	 -0.00000273935060780554  -0.00000020574799465699	      -0.00000004363715847430	  	-0.00000001430608195426 	-0.00000000591587973546 &  
-1.588071022611375312718684509421    -0.004602881871269254502543511851 -0.000002739351068678330086744127    -0.000000205753231647526892605590 	-0.000000043659496185302904545860  -0.000000014366823059961810214623  -0.000000006034363442475256753121 & 
-1.588071022611375312718684509423  -0.004602881871269254502543511851   -0.000002739351068678330086744132  -0.000000205753231647526892605535  -0.000000043659496185302904545958  -0.000000014366823059961810213641 -0.000000006034363442475256759610 \\
\hline
\end{tabular}
\end{table}
\end{center}
%%%%%%%%%%%%%%%%%%%%%%%%%%%%%%%%%%%%%%%%%%

\end{document}